\documentclass[11pt]{article}
\usepackage{amscd}
\usepackage{amsfonts}
\usepackage{amsmath}
\usepackage{amssymb}
\usepackage{amsthm}
\usepackage{bbm}
\usepackage{CJK}
\usepackage{fancyhdr}
\usepackage{graphicx}
\usepackage{hyperref}
\usepackage{indentfirst}
\usepackage{latexsym}
\usepackage{mathrsfs}
\usepackage{xypic}

\newtheorem{theorem}{Theorem}[section]
\newtheorem{lemma}[theorem]{Lemma}
\newtheorem{definition}[theorem]{Definition}
\newtheorem{proposition}[theorem]{Proposition}

\newtheorem{remark}[theorem]{Remark}

\usepackage[top=1in,bottom=1in,left=1.25in,right=1.25in]{geometry}
\def\<{\langle}
\def\>{\rangle}
\def\a{\alpha}
\def\b{\beta}

\def\c{\cdot}

\date{}
\begin{document}
\renewcommand{\baselinestretch}{1.2}
\renewcommand{\arraystretch}{1.0}
\title{\bf Representations and deformations of Hom-Lie-Yamaguti superalgebras}
\author{{\bf Shuangjian Guo$^{1}$, Xiaohui Zhang$^{2}$,  Shengxiang Wang$^{3}$\footnote
        { Corresponding author(Shengxiang Wang):~~wangsx-math@163.com} }\\
{\small 1. School of Mathematics and Statistics, Guizhou University of Finance and Economics} \\
{\small  Guiyang  550025, P. R. of China} \\
{\small 2.  School of Mathematical Sciences, Qufu Normal University}\\
{\small Qufu  273165, P. R. of China}\\
{\small 3.~ School of Mathematics and Finance, Chuzhou University}\\
 {\small   Chuzhou 239000,  P. R. of China}}
 \maketitle
\begin{center}
\begin{minipage}{13.cm}

{\bf \begin{center} ABSTRACT \end{center}}
Let $(L, \a)$  be a Hom-Lie-Yamaguti superalgebra.  We first introduce  the representation and cohomology theory of Hom-Lie-Yamaguti superalgebras. Furthermore, we introduce the notions of generalized derivations and representations of $(L, \a)$  and present  some properties. Finally, we  investigate the deformations of $(L, \a)$  by choosing some suitable cohomology.
 \smallskip

{\bf Key words}: Hom-Lie-Yamaguti superalgebra; representation; cohomology;   derivation;  deformation.
 \smallskip

 {\bf 2010 MSC:} 17A40; 17A30, 17B10; 17B56
 \end{minipage}
 \end{center}
 \normalsize\vskip0.5cm

\section{Introduction}
\def\theequation{\arabic{section}. \arabic{equation}}
\setcounter{equation} {0}

Lie triple systems arosed initially in Cartan's study of Riemannian geometry.
Jacobson \cite{Jacobson} first introduced Lie triple systems and Jordan triple systems in connection with problems from Jordan theory
and quantum mechanics, viewing Lie triple systems as subspaces of Lie algebras that are closed relative to the ternary product.
Lie-Yamaguti algebras were introduced by Yamaguti in \cite{Yamaguti1958} to give an algebraic interpretation of the characteristic properties of the torsion and curvature of homogeneous
spaces with canonical connection in \cite{Nomizu}. He called them ¡°generalized Lie triple systems¡±
at first, which were later called ``Lie triple algebras". Recently, they were renamed as
``Lie-Yamaguti algebras".
 \medskip

In \cite{Gaparayi2012}, the authors introduced the concept of Hom-Lie-Yamaguti algebras. It is a Hom-type
generalization of a Lie-Yamaguti algebra.  In \cite{Ma}, the authors studied the formal deformations of Hom-Lie-Yamaguti
algebras.  In \cite{Zhang2015}, the authors introduced the representation and cohomology theory of Hom-Lie-Yamaguti algebras, and studied deformation and extension of Hom-Lie-Yamaguti algebras an an application.  In \cite{Mandal}, the authors introduced the notion of extensions of Hom-Lie-Rinehart algebras and deduced a characterisation of low dimensional
cohomology spaces in terms of the group of automorphisms of certain abelian extension and the equivalence classes of those abelian extensions in the category of Hom-Lie-Rinehart algebras.  In \cite{Zhang2018}, the authors introduced the notion of crossed modules for Hom-Lie-Rinehart algebras, and   studied  their construction  of Hom-Lie-Rinehart algebras.

In \cite{Gaparayi2019}, the authors introduced   the concept of Hom-Lie-Yamaguti superalgebras  and give some examples of Hom-
Lie-Yamaguti superalgebras. Later, in \cite{Gaparayi2018}, the authors studied the relation of Hom-Leibniz superalgebras and  Hom-Lie-Yamaguti superalgebras.

 The purpose of this paper is to extend  the Hom-type generalization of binary superalgebras to the one of
ternary superalgebras or binary-ternary superalgebras.  This paper is organized as follows.
In Section 2, we recall the definition of Hom-Lie-Yamaguti superalgebras.
In Section 3, We  introduce  the representation and cohomology theory of Hom-Lie-Yamaguti superalgebras.
In Section 4 we introduce the notions of generalized derivations and representations of a Hom-Lie-Yamaguti superalgebra  and present  some properties.
In Section 5, we consider  the theory of deformations of a Hom-Lie-Yamaguti superalgebra  by choosing a suitable cohomology.

\section{Preliminaries}
\def\theequation{\arabic{section}.\arabic{equation}}
\setcounter{equation} {0}

Throughout this paper, we work on an algebraically closed field $\mathbb{K}$ of characteristic different
from 2 and 3. We  recall some basic definitions regarding Hom-Lie-Yamaguti superalgebras from \cite{Gaparayi2019}.

\begin{definition}
A Hom-Lie-Yamaguti superalgebra (Hom-LY superalgebra for short) is a quadruple $(L, [\c, \c], \{\c, \c, \c\}, \a)$ in which $L$ is $\mathbb{K}$-vector superspace, $[\c, \c]$ a binary superoperation and $\{\c, \c, \c\}$
a ternary superoperation on $L$, and $\a: L\rightarrow L$ an even linear map such that
\begin{eqnarray*}
&&(SHLY1)  ~~\a([x, y])=[\a(x), \a(y)],\\
&&(SHLY2) ~~\a(\{x,y,z\})=\{\a(x),\a(y),\a(z)\},\\
&&(SHLY3)~~ [x, y]=-(-1)^{|x||y|}[y,  x],\\
&&(SHLY4) ~~\{x,y,z\}=-(-1)^{|x||y|}\{y,x,z\},\\
&&(SHLY5) ~~\circlearrowleft_{(x,y,z)}(-1)^{|x||z|}([[x, y], \a(z)]+\{x, y, z\})=0,\\
&&(SHLY6) ~~\circlearrowleft_{(x,y,z)}(-1)^{|x||z|}(\{[x, y], \a(z), \a(u)\})=0,\\
&&(SHLY7) ~~\{\a(x), \a(y), [u, v]\}=[\{x, y, u\}, \a^{2}(v)]+(-1)^{|u|(|x|+|y|)}[\a^{2}(u),\{x, y, v\}],\\
&&(SHLY8) ~~ \{\a^2(x), \a^{2}(y), \{u, v, w\}\}=\{\{x, y, u\}, \a^{2}(v), \a^2(w)\}\\
&&\hspace{6.5cm} +(-1)^{|u|(|x|+|y|)}\{\a^{2}(u),  \{x, y, v\}, \a^{2}(w)\}\\
&&\hspace{6.5cm} +(-1)^{(|u|+|v|)(|x|+|y|)}\{\a^{2}(u),  \a^{2}(v), \{x, y, w\}\},
\end{eqnarray*}
for all $ x, y, z, u, v, w\in L$ and where $\circlearrowleft_{(x,y,z)}$ denotes the sum over cyclic permutation of $x, y, z$, and $|x|$ denotes the degree of the element $x\in L$. We denote a Hom-LY superalgebra by $(L, \a)$.
\end{definition}

A homomorphism between two Hom-LY superalgebras  $(L, \a)$ and $(L', \a')$ is a linear map $\varphi: L\rightarrow L'$ satisfying $\varphi \circ \a=\a\circ \varphi$ and
\begin{eqnarray*}
\varphi([x, y])=[\varphi(x), \varphi(y)]', ~~~~~~~~~\varphi(\{x, y, z\})=\{\varphi(x), \varphi(y),\varphi(z)\}'.
\end{eqnarray*}

\begin{remark} (1) If $\a=Id$, then the Hom-LY superalgebra $(L, [\c, \c], \{ \c, \c, \c \}, \a)$ reduces to a LY superalgebra $(L, [\c, \c], \{\c, \c, \c\})$ (see (SLY 1)-(SLY 6)).

(2) If $[x, y ]= 0$, for all $x, y \in L$, then $(L, [\c, \c], \{\c, \c, \c\}, \a)$ becomes a Hom-Lie supertriple system
$(L,  \{\c, \c, \c\}, \a^{2})$.

(3) If $\{x, y, z\} =0$ for all $x, y, z\in L$, then the Hom-LY superalgebra $(L, [\c, \c], \{\c, \c, \c\}, \a)$ becomes a
Hom-Lie superalgebra $(L, [\c, \c], \a)$.

\end{remark}

\section{Representations  of Hom-Lie-Yamaguti superalgebras}
\def\theequation{\arabic{section}. \arabic{equation}}
\setcounter{equation} {0}

\begin{definition} Let $(L, \a)$ be a Hom-LY superalgebra and $(V, \b)$ be a Hom-vector space. A representation of
$L$ on $V$ consists of a linear map $\rho: L\rightarrow$ End($V$) and bilinear maps $D, \theta: L\times L\rightarrow$ End($V$)
such that the following conditions are satisfied:
\begin{eqnarray*}
 &&(SHR1) ~\rho(\a(x))\circ \beta =\beta\circ \rho(\a(x)),\\
 &&(SHR2) ~D(\a(x), \a(y))\circ \b=\b\circ D(x, y),\\
 &&(SHR3) ~\theta(\a(x), \a(y))\circ \b=\b\circ \theta(x, y),\\
 &&(SHR4) ~D(x, y)-(-1)^{|x||y|}\theta(y, x)+\theta(x, y)+\rho([x, y])\circ \beta\\
  &&\hspace{6cm}-\rho(\a(x))\rho(y)+(-1)^{|x||y|}\rho(\a(y))\rho(x)=0,\\
 &&(SHR5) ~D([x, y], \a(z))+(-1)^{|x|(|y|+|z|)}D([y, z], \a(x))+(-1)^{|z|(|x|+|y|)}D([z, x], \a(y))=0, \\
 &&(SHR6)  ~\theta([x, y], \a(z))\circ \b=(-1)^{|y||z|}\theta(\a(x), \a(z))\rho(y)-(-1)^{|x|(|y|+|z|)}\theta(\a(y), \a(z))\rho(x), \\
 &&(SHR7) ~D(\a(x), \a(y))\rho(z)=(-1)^{|z|(|x|+|y|)}\rho(\a^{2}(z))D(x, y)+\rho(\{x, y, z\})\circ \b^{2}, \\
 &&(SHR8)~\theta(\a(x), [y, z])\circ \b=(-1)^{|x||y|}\rho(\a^{2}(y))\theta(x, z)-(-1)^{|z|(|x|+|y|)}\rho(\a^{2}(z))\theta(x, y),\\
 &&(SHR9) ~D(\a^{2}(x), \a^{2}(y))\theta(u, v)=~(-1)^{(|u|+|v|)(|x|+|y|)}\theta(\a^{2}(u), \a^{2}(v))D(x, y)\\
 &&\hspace{4cm}+\theta(\{x, y, u\}, \a^{2}(v))\circ \b^{2}+(-1)^{|u|(|x|+|y|)}\theta(\a^{2}(u), \{x, y, v\})\circ \b^{2},\\
 &&(SHR10)  ~\theta(\a^{2}(x), \{y, z, u\})\circ \b^{2}=(-1)^{(|z|+|u|)(|x|+|y|)}\theta(\a^{2}(z), \a^{2}(u)) \theta(x, y)\\
 &&\hspace{2cm}-(-1)^{|y||z|}\theta(\a^{2}(y), \a^{2}(u)) \theta(x, z)+(-1)^{|x|(|y|+|z|)}D(\a^{2}(y), \a^{2}(z))\theta(x, u),
\end{eqnarray*}
for any $x, y, z, u, v\in L$. In this case, we also call $V$ to be a $L$-module.

\end{definition}
\begin{proposition}
Let $(L, \a)$ be a Hom-LY superalgebra and $V$ be a $\mathbb{Z}_2$-graded vector space. Assume we have a
map $\rho$ from $L$ to End($V$) and maps $D, \theta: L\times L\rightarrow$ End($V$). Then $(\rho, D, \theta)$ is a representation
of $(L, \a)$ on $(V, \b)$ if and only if $L \oplus V$ is a Hom-LY superalgebra under the following maps:
\begin{eqnarray*}
&&(\a+\b)(x+u):=\a(x)+\b(u),\\
&&[x+u, y+v]:=[x, y]+\rho(x)(v)- (-1)^{|x||y|}\rho(y)(u),\\
&&\{x+u, y+v, z+w\}:=\{x, y, z\}+D(x, y)(w)\\
&&\hspace{4cm}-(-1)^{|y||z|}\theta(x, z)(v)+(-1)^{|x|(|y|+|z|)}\theta(y, z)(u),
\end{eqnarray*}
for any $x, y, z\in L$ and $u, v, w\in V$.
\end{proposition}
{\bf Proof.} It is easy to check that the conditions (SHLY1)-(SHLY4) hold, we only verify that conditions (SHLY5)-(SHLY8) hold for maps defined on $L \oplus V$.

For (SHLY5), we have
\begin{eqnarray*}
&&\{x+u, y+v, z+w\}+c.p.\\
&=&[\{x, y, z\}+D(x, y)(w)\\
&&\hspace{3cm}-(-1)^{|y||z|}\theta(x, z)(v)+(-1)^{|x|(|y|+|z|)}\theta(y, z)(u)]+c.p.
\end{eqnarray*}
and
\begin{eqnarray*}
&&[[x+u, y+v], \a(z)+\b(w)]+c.p.\\
&=&[[x, y]+\rho(x)(v)- (-1)^{|x||y|}\rho(y)(u), \a(z)+\b(w)]+c.p.\\
&=&([[x, y], \a(z)]+\rho([x, y])\b(w) -(-1)^{|x||z|}\rho(\a(z))\rho(x)(v)\\
&&\hspace{6cm}+(-1)^{(|x|+|z|)|y|}\rho(\a(z))\rho(y)(u))+c.p.
\end{eqnarray*}
Thus by (SHR4), the condition (SHLY5) holds.

For (SHLY6), we have
\begin{eqnarray*}
&&(-1)^{|x||z|}(\{[x+u, y+v], \a(z)+\b(w), \a(p)+\b(t)\})+c.p.\\
&=& (-1)^{|x||z|}(\{[x, y], \a(z), \a(p)\}+D([x, y], \a(z))(\b(t))-\theta([x, y], \a(p))(\b(w))\\
&&\hspace{6cm}+(-1)^{(|z|+|p|+|y|)|x|}\theta(\a(z), \a(p))(\rho(x)(v))\\
&&\hspace{6cm}-(-1)^{(|z|+|p|+|y|)|x|}\theta(\a(z), \a(p))(\rho(y)(u)))+c.p.\\
&=& 0.
\end{eqnarray*}
Thus by (SHR5), the condition (SHLY6) holds.

For (SHLY7), we have
\begin{eqnarray*}
&&\{\a(x)+\b(u), \a(y)+\b(v), [z+w, p+t]\}\\
&=& \{\a(x), \a(y), [z, p]\}+D(\a(x), \a(y))(\rho(z)(t))- (-1)^{|z||p|}D(\a(x), \a(y))(\rho(p)(w))\\
&& -(-1)^{(|z|+|p|)|y|}\theta(x, [z, p])(v)+(-1)^{(|z|+|p|+|y|)|x|}\theta(y, [z, p])(u),
\end{eqnarray*}
and
\begin{eqnarray*}
&& [\{x+u, y+v, z+w\},  \a^{2}(p)+\b^{2}(t)]\\
&&\hspace{5cm}+(-1)^{|z|(|x|+|y|)}[ \a^{2}(z)+\b^{2}(w), \{x+u, y+v, p+t\}]\\
&=&[\{x, y, z\}, \a^{2}(p)]+\rho(\{x, y, z\})(t)-(-1)^{(|x|+|y|+|z|)|p|}\rho(p)(D(x, y)(w))\\
&&-(-1)^{(|x|+|y|+|z|)|p|+|y||z|}\rho(p)(\theta(x, z)(v))+(-1)^{(|z|+|y|)(|p|+|x|)+|p||x|}\rho(p)(\theta(y, z)(u))\\
&&+(-1)^{|z|(|x|+|y|)}[\a^{2}(z), \{x, y, p\}]+(-1)^{|z|(|x|+|y|)}\rho(z)(D(x, y)(t))\\
&&-(-1)^{|x||z|+(|x|+|p|)|y|}\rho(z)(\theta(x, p)(v))+(-1)^{|z||p|+|x|(|z|+|y|+|p|)}\rho(z)(\theta(y, p)(u))\\
&&\hspace{9cm}-(-1)^{|p||z|}\rho(\{x, y, p\})(w).
\end{eqnarray*}
Thus by (SHR6), the condition (SHLY7) holds.

Now it suffices to verify (SHLY8). By the definition of a Hom-LY superalgebra, we have
\begin{eqnarray*}
&&\{\a^2(x_1)+\b^{2}(u_1), \a^{2}(x_2)+\b^{2}(u_2), \{y_1+v_1, y_2+v_2, y_3+v_3\}\}\\
&=& \{\a^2(x_1), \a^{2}(x_2), \{y_1, y_2, y_3\}\}-(-1)^{|x_2|(|y_1|+|y_2|+|y_3|)}\theta(\a^2(x_1), \{y_1, y_2, y_3\})(\b^{2}(u_2))\\
&&+(-1)^{|x_1|(|x_2|+|y_1|+|y_2|+|y_3|)}\theta(\a^2(x_2), \{y_1, y_2, y_3\})(\b^{2}(u_1))\\
&&+ D(\a^2(x_1), \a^{2}(x_2))(D(y_1,y_2)(v_3))-(-1)^{|y_2||y_3|}D(\a^2(x_1), \a^{2}(x_2))(\theta(y_1, y_3)(v_2))\\
&&+(-1)^{|y_1|(|y_2|+|y_3|)}D(\a^2(x_1), \a^{2}(x_2))(\theta(y_2, y_3)(v_1)),\\
&&\{\{x_1+u_1, x_2+u_2, y_1+v_1\}, \a^2(y_2)+\b^{2}(v_2), \a^{2}(y_3)+\b^{2}(v_3)\}\\
&=& \{\{x_1,x_2,y_1\}, \a^2(y_2)+\a^{2}(y_3)\}+D(\{x_1,x_2,y_1\}, \a^{2}(y_2))(\b^{2}(v_3))\\
&& +(-1)^{|y_2||y_3|}\theta(\{x_1,x_2,y_1\}, \a^{2}(y_3))(\b^{2}(u_1))+(-1)^{(|y_2|+|y_3|)(|x_1|+|x_2|+|y_1|)}\theta( \a^2(y_2), \a^{2}(y_3))\\
&&(D(x_1,x_2)(v_1))-(-1)^{(|y_1|+|y_2|+|y_3|)|x_2|+(|x_1|+|y_1|)(|y_2|+|y_3|)}\theta( \a^2(y_2), \a^{2}(y_3))(\theta(x_1,y_1)(u_2))\\
&&+(-1)^{(|y_2|+|y_3|)(|x_2|+|y_3|)+|x_1|(|y_2|+|y_3|+|x_2|+|y_3|)}\theta( \a^2(y_2), \a^{2}(y_3))(\theta(x_2,y_1)(u_1)),\\
&&(-1)^{|y_1|(|x_1|+|x_2|)}\{\a^2(y_1)+\b^{2}(v_1), \{x_1+u_1, x_2+u_2, y_2+v_2\}, \a^{2}(y_3)+\b^{2}(v_3)\}\\
&=& (-1)^{|y_1|(|x_1|+|x_2|)}\{\a^2(y_1), \{x_1, x_2, y_2\}, \a^{2}(y_3)\}+(-1)^{|y_1|(|x_1|+|x_2|)}D(\a^2(y_1), \{x_1,x_2, y_2\})(\b^{2}(v_3))\\
&& +(-1)^{|y_1|(|y_2|+|y_3|)}\theta(\{x_1,x_2, y_2\}, \a^2(y_3))(\b^{2}(v_1))-(-1)^{|y_2||y_3|}\theta( \a^2(y_1), \a^{2}(y_3))(D(x_1,x_2)(v_3))\\
&&+(-1)^{(|y_1|+|y_3|)(|x_1|+|y_2|)+|x_2|(|y_1|+|y_2|+|x_3|)+|y_1||y_2|}\theta( \a^2(y_1), \a^{2}(y_3))(\theta(x_1,y_2)(u_2))\\
&&-\theta( \a^2(y_1), \a^{2}(y_3))(\theta(x_2,y_2)(u_1))\\
&&(-1)^{(|y_1|+|y_2|)(|x_1|+|x_2|)}\{\a^2(y_1)+\b^{2}(v_1),  \a^{2}(y_2)+\b^{2}(v_2), \{x_1+u_1, x_2+u_2, y_3+v_3\}\}\\
&=& (-1)^{(|y_1|+|y_2|)(|x_1|+|x_2|)}[\a^2(y_1), \a^{2}(y_2), \{x_1,x_2,y_3\}]-(-1)^{|y_1|(|x_1|+|x_2|)+|y_2||y_3|}\\
&&-\theta(\a^{2}(y_1), \{x_1,x_2,y_3\})(\b^{2}(v_2))+(-1)^{|y_2|(|x_1|+|x_2|)+|y_1|(|y_2|+|y_3)|}\theta(\a^{2}(y_2), \{x_1,x_2,y_3\})(\b^{2}(v_1))\\
&& +(-1)^{(|y_1|+|y_2|)(|x_1|+|x_2|)}D( \a^2(y_1), \a^{2}(y_2))(D(x_1,x_2)(v_3))\\
&&-(-1)^{(|y_1|+|y_2|)(|x_1|+|x_3|)+|x_2||y_3|}D( \a^2(y_1), \a^{2}(y_2))(\theta(x_1,y_3)(u_2))\\
&&+(-1)^{(|y_1|+|y_2|)(|x_2|+|y_3|)+|x_1|(|x_2|+|y_1|+|y_2|+|y_3|)}D( \a^2(y_1), \a^{2}(y_2))(\theta(x_2,y_3)(u_1))
\end{eqnarray*}
Thus by (SHR7), the condition (SHLY8) holds.  Therefore we obtain that a Hom-LY superalgebra on $L\oplus V$.  $\hfill \Box$

Let $V$ be a representation of Hom-LY superalgebra $L$. Let us define the cohomology groups of $L$ with coefficients
in $V$. Let $f: L\times L \times \c\c\c\times L$ be $n$-linear maps of $L$ into $V$ such that the following conditions are
satisfied:
\begin{eqnarray*}
&&f(\a(x_1), \c\c\c, \a(x_n))=\beta(f(x_1, \c \c \c, x_n)),\\
&& f(x_1, \c \c \c, x_{2i-1}, x_{2i}, \c \c \c, x_{n})=0, ~~~~\mbox{if}~~~ x_{2i-1}=x_{2i}.
\end{eqnarray*}
The vector space spanned by such linear maps is called an $n$-cochain of $L$, which is denoted by
$C^{n}(L, V)$ for $n\geq 1$.

\begin{definition}
For any $(f, g)\in  C^{n}(L, V) \times C^{2n+1}(L, V )$ the coboundary operator  $\delta: (f, g)\rightarrow (\delta_I f, \delta_{II} g)$ is a mapping from $C^{2n}(L, V) \times C^{2n+1}(L, V )$ into $C^{2n+2}(L, V) \times C^{2n+3}(L, V )$ defined as follows:
\begin{eqnarray*}
&& (\delta_If)(x_1,x_2,\c\c\c,x_{2n+2})\\
&=&\\
&& (\delta_{II}g)(x_1,x_2,\c\c\c,x_{2n+3})\\
&=& 
\end{eqnarray*}
\end{definition}

\begin{proposition}

The coboundary operator defined above satisfies $\delta \circ \delta=0$, that is $\delta_I\circ \delta_I=0$ and $\delta_{II}\circ \delta_{II}=0$.

\end{proposition}
{\bf Proof.} Similar to \cite{Zhang2015}. $\hfill \Box$

Let $Z^{2n}(L, V) \times Z^{2n+1}(L, V )$ be the subspace of $C^{2n}(L, V) \times C^{2n+1}(L, V )$ spanned by $(f, g)$
such that $\delta(f, g) = 0$ which is called the space of cocycles and $B^{2n}(L, V ) \times B^{2n+1}(L, V ) =
\delta(C^{2n-2}(L, V ) \times C^{2n-1}(L, V ))$ which is called the space of coboundaries.

\begin{definition}
For the case $n\geq 2$, the $(2n, 2n + 1)$-cohomology group of a Hom-LY superalgebra $L$ with
coefficients in $V$ is defined to be the quotient space:
\begin{eqnarray*}
 H^{2n}(L, V)\times H^{2n+1}(L, V):=(Z^{2n}(L, V) \times Z^{2n+1}(L, V ))/(B^{2n}(L, V) \times B^{2n+1}(L, V )).
\end{eqnarray*}
\end{definition}
In conclusion, we obtain a cochain complex whose cohomology group is called cohomology
group of a Hom-LY superalgebra $L$ with coefficients in $V$.
\section{Derivations  of Hom-Lie-Yamaguti superalgebras}
\def\theequation{\arabic{section}. \arabic{equation}}
\setcounter{equation} {0}

In this section, we give the definition of  derivations of  Hom-LY superalgebras,
then we study its generalized derivations.

\begin{definition}
 A  linear map $D: L\rightarrow L$  is called the  $\a^k$-derivation of $L$  if it satisfies
\begin{eqnarray*}
&&D([x, y])=(-1)^{|D||x|}[\a^{k}(x), D(y)]+[D(x), \a^{k}(y)],\\
&&D(\{x,y,z\})=\{D(x),\a^{k}(y),\a^{k}(z)\}+(-1)^{|D||x|}\{\a^{k}(x), D(y), \a^{k}(z)\}\nonumber\\
&&\hspace{6cm}+(-1)^{|D|(|x|+|y|)}\{\a^{k}(x),\a^{k}(y),D(z)\},~~~~~~~
\end{eqnarray*}
for all $x,y,z\in L$, where $|D|$ denotes the degree of $D$.
\end{definition}

We denote by $Der(L)=\bigoplus_{k\geq 0}Der_{\a^k}(L)$, where $Der_{\a^k}(L)$ is the set of all homogeneous $\a^k$-derivations of $L$.
Obviously, $Der(L)$ is a subalgebra of $End(L)$ and has a normal Lie superalgebra structure via the bracket product
\begin{eqnarray*}
[D,D']=DD'-(-1)^{|D||D'|}D'D.
\end{eqnarray*}

\begin{theorem}
Set $Der(L)=\bigoplus_{k\geq 0}Der_{\a^{k}}(L)$. Then $Der(L)$ is a Lie superalgebra.

\end{theorem}
{\bf Proof.} It is sufficient to prove $[Der_{\a^k}(L), Der_{\a^s}(L)]\subseteq Der_{\a^{k+s}}(L)$. It is easy to check that $[D, D']\circ \a=\a\circ [D, D']$.

Note that
\begin{eqnarray*}
&&[D, D']([x, y])\\
&=&D([\a^s(x), D'(y))]+(-1)^{|D'||x|}[D'(x), \a^{s}(y)])\\
&&\hspace{5cm}-(-1)^{|D||D'|}D'([\a^{k}(x), D(y)]+(-1)^{|D||x|}[D(x), \a^{k}(y)])\\
&=&(-1)^{|D'||x|}[D\a^{s}(x), \a^{k}D'(y)]+(-1)^{|D||x|+|D'||x|}[\a^{k+s}(x), DD'(y)]+[DD'(x), \a^{k+s}(y)]\\
&&+(-1)^{|D|(|D'|+|x|)}[\a^{k}D'(x), D\a^{s}(y)] -(-1)^{|D||x|+|D||D'|}[D'\a^k(x), \a^sD(y)]\\
&&-(-1)^{|D'||x|+|D||x|+|D||D'|}[\a^{k+s}(x), D'D(y)]-(-1)^{|D||D'|}[D'D(x), \a^{k+s}(y)]\\
&&\hspace{6cm}-(-1)^{|D'|(|x|+|D|)+|D||D'|}[\a^sD(x), D'\a^{k}(y)]\\
&=& [[D, D'](x), \a^{k+s}(y)]+(-1)^{(|D'|+|D|)|x|}[\a^{k+s}(x), [D,D'](y)].
\end{eqnarray*}
Similarly, we can check that
\begin{eqnarray*}
&&[D, D'](\{x, y, z\})\\
&=& \{[D, D'](x), \a^{k+s}(y), \a^{k+s}(z)\}+(-1)^{(|D'|+|D|)|x|}\{\a^{k+s}(x), [D, D'](y),  \a^{k+s}(z)\}\\
&&\hspace{4cm}+(-1)^{(|D'|+|D|)(|x|+|y|)}\{\a^{k+s}(x),  \a^{k+s}(y),  [D, D'](z)\}.
\end{eqnarray*}
It follows that $[D, D']\in Der_{\a^{k+s}}(L)$. $\hfill \Box$
\begin{definition}
Let $(L, \a)$  be a Hom-Lie-Yamaguti superalgebra.
$D\in End_{\overline{s}}(L)$ is said to be a homogeneous generalized $\a^k$-derivation of $L$, if there  exist three endomorphisms
$D', D'',D'''\in End_{\overline{s}}(L)$ such that
\begin{eqnarray*}
&&[D(x), \a^{k}(y)]+(-1)^{s|x|}[\a^{k}(x), D'(y)]=D''([x, y]),\\
&&\{D(x),\a^{k}(y),\a^{k}(z)\}+(-1)^{s|x|}\{\a^{k}(x), D'(y),\a^{k}(z)\}\\
&&\hspace{4cm}+(-1)^{s(|x|+|y|)}\{\a^{k}(x),\a^{k}(y),D''(z)\}=D'''(\{x,y,z\}),
\end{eqnarray*}
for all $x,y,z\in L$.
\end{definition}
\begin{definition}
Let $(L, \a)$  be a Hom-Lie-Yamaguti superalgebra.
$D\in End_{\overline{s}}(L)$ is said to be a homogeneous  $\a^k$-quasiderivation of $L$, if there  exist an endomorphism
$D', D''\in End_{\overline{s}}(L)$ such that
\begin{eqnarray*}
&&[D(x), \a^{k}(y)]+(-1)^{s|x|}[\a^{k}(x), D(y)]=D'([x, y]),\\
&&\{D(x),\a^{k}(y),\a^{k}(z)\}+(-1)^{s|x|}\{\a^{k}(x), D(y),\a^{k}(z)\}\\
&&\hspace{4cm}+(-1)^{s(|x|+|y|)}\{\a^{k}(x),\a^{k}(y),D(z)\}=D''(\{x,y,z\}),
\end{eqnarray*}
for all $x,y,z\in L$.
\end{definition}

Let $GDer(L)$ and $QDer(L)$ be the sets of homogeneous generalized $\a^k$-derivations and of homogeneous  $\a^k$-quasiderivations, respectively. That is,
\begin{eqnarray*}
GDer(L):=\bigoplus_{k\geq 0}GDer_{\a^k}(L),~QDer(L):=\bigoplus_{k\geq 0}QDer_{\a^k}(L).
\end{eqnarray*}

\begin{definition}
Let $(L, \a)$  be a Hom-Lie-Yamaguti superalgebra.
The $\a^k$-centroid of  $L$ is the space of linear transformations on $L$ given by
\begin{eqnarray*}
&&C_{\a^k}(L)=\{D\in End(L)| [D(x), \a^{k}(y)]=(-1)^{|D||x|}[\a^{k}(x), D(y)]=D([x, y])~~\mbox{and}~~  \\
&& \{D(x),\a^{k}(y),\a^{k}(z)\}=(-1)^{|D||x|}\{\a^{k}(x),D(y),\a^{k}(z)\}\\
&&\hspace{3.5cm}=(-1)^{|D|(|x|+|y|)}\{\a^{k}(x),\a^{k}(y),D(c)\}=D(\{x,y,z\})\}.
\end{eqnarray*}
We denote $C(L)=\bigoplus_{k\geq 0}C_{\a^k}(L)$  and call it the centroid of  $L$.
\end{definition}

\begin{definition}
Let $(L, \a)$  be a Hom-Lie-Yamaguti superalgebra.
The  quasicentroid of $L$ is the space of linear transformations on $L$ given by
\begin{eqnarray*}
&&QC_{\a^{k}}(L)=\{D\in End(L)|[D(x), \a^{k}(y)]=(-1)^{|D||x|}[\a^{k}(x), D(y)]~~\mbox{and}~~     \{D(x),\a^{k}(y),\a^{k}(z)\}\\
&&=(-1)^{|D||x|}\{\a^{k}(x),D(y),\a^{k}(z)\}=(-1)^{|D|(|x|+|y|)}\{\a^{k}(x),\a^{k}(y),D(c)\}\},
\end{eqnarray*}
for all $x,y,z\in L$. We denote $QC(L)=\bigoplus_{k\geq 0}QC_{\a^k}(L)$  and call it the quasicentroid  of  $L$.
\end{definition}

\begin{remark}
Let $(L, \a)$  be a Hom-Lie-Yamaguti superalgebra.
Then $C(L)\subseteq QC(L).$
\end{remark}

\begin{definition}
Let $(L, \a)$  be a Hom-Lie-Yamaguti superalgebra.
$D\in End(L)$ is said to be a central derivation  of $L$ if
\begin{eqnarray*}
&&[D(x), \a^{k}(y)]=D([x, y])=0,\\
&&D([x,y,z])=[D(x),\a^{k}(y),\a^{k}(z)]=0,
\end{eqnarray*}
for all $x,y,z\in L$.
Denote the set of all central  derivations by $ZDer(L):=\bigoplus_{k\geq 0}Der_{\a^k}(L)$.
\end{definition}

\begin{remark}
Let $(L, \a)$  be a Hom-Lie-Yamaguti superalgebra.
Then
\begin{eqnarray*}
ZDer(L)\subseteq Der(L)\subseteq QDer(L)\subseteq GDer(L) \subseteq End(L).
\end{eqnarray*}
\end{remark}
\begin{definition}
Let $(L, \a)$  be a Hom-Lie-Yamaguti superalgebra.
If $Z(L)=\{x\in L|~[x,y,z]=0,~\forall~x, y, z\in L\}$,
then  $Z(L)$ is called the center of $L$.
\end{definition}

\begin{proposition}
Let $(L, \a)$  be a Hom-Lie-Yamaguti superalgebra,
 then the following statements hold:

(1)~$GDer(L), QDer(L)$ and $C(L)$ are subalgebras of $End(L)$.

(2)~$ZDer(L)$ is an ideal of $Der(L)$.
\end{proposition}
{\bf Proof.}
(1) We only prove that $GDer(L)$ is a subalgebra of $End(L)$, and similarly for cases of $ QDer(L)$ and $C(L)$.
For any $D_1\in GDer_{\a^k}(L),D_2\in GDer_{\a^s}(L)$ and $x,y,z\in L$, we have
\begin{eqnarray*}
&&\{D_1D_2(x),\a^{k+s}(y),\a^{k+s}(z)\}\\
&=&D'''_1\{D_2(x),\a^{s}(y),\a^{s}(z)\}-(-1)^{|D_1|(|D_2|+|x|)}\{\a^{k}(D_2(x)),D'_1(\a^{s}(y)),\a^{k+s}(z)\}\\
  &&-(-1)^{|D_1|(|D_2|+|x|+|y|)}\{\a^{k}(D_2(x)),\a^{k+s}(y),D''_1(\a^{s}(z))\}\\
&=& D'''_1\{ D'''_2(\{x,y,z\})-(-1)^{|D_2||x|}\{\a^{s}(x), D'_2(y),\a^{s}(z)\}-(-1)^{|D_2|(|x|+|y|)}\{\a^{s}(x),\a^{s}(y),D''_2(z)\}\}\\
  &&-(-1)^{|D_1|(|D_2|+|x|)}\{\a^{k}(D_2(x)),D'_1(\a^{s}(y)),\a^{k+s}(z)\}\\
  &&-(-1)^{|D_1|(|D_2|+|x|+|y|)}\{\a^{k}(D_2(x)),\a^{k+s}(y),D''_1(\a^{s}(z))\}\\
&=&D'''_1D'''_2(\{x,y,z\})-(-1)^{|D_2||x|}D'''_1\{x, D'_2(y),c\}-(-1)^{|D_2|(|x|+|y|)}D'''_1\{x,y,D''_2(z)\}\\
  &&-(-1)^{|D_1|(|D_2|+|x|)}\{D_2(x),D'_1(y),z\}-(-1)^{|D_1|(|D_2|+|x|+|y|)}\{D_2(x),y,D''_1(z)\}\\
&=&D'''_1D'''_2(\{x,y,z\})
   -(-1)^{|D_2||x|}\{D_1(\a^{s}(x)), \a^{k}(D'_2(y)),\a^{k+s}(z)\}\\
   &&-(-1)^{(|D_1|+|D_2|)|x|}\{\a^{k+s}(x), D'_1 D'_2(y),\a^{k+s}(z)\}\\
   &&-(-1)^{|D_2||x|+|D_1|(|x|+|y|+|D_2|)}\{\a^{k+s}(x), \a^{k}(D'_2(y)),D''_1(\a^{s}(z))\}\\
   &&-(-1)^{|D_2|(|x|+|y|)}\{D_1(\a^{s}(x)),\a^{k+s}(y),\a^{k}(D''_2(z))\}\\
   &&-(-1)^{|D_2|(|x|+|y|)+|D_1||x|}\{\a^{k+s}(x),D'_1(\a^{s}(y)),\a^{k}(D''_2(z))\}\\
   &&-(-1)^{(|D_1|+|D_2|)(|x|+|y|)}\{\a^{k+s}(x),\a^{k+s}(y),D''_1 D''_2(z)\}\\
   &&-(-1)^{|D_1|(|D_2|+|x|)}\{\a^{k}(D_2(x)),D'_1(\a^{s}(y)),\a^{k+s}(z)\}\\
   &&-(-1)^{|D_1|(|D_2|+|x|+|y|)}\{\a^{k}(D_2(x)),\a^{k+s}(y),D''_1(\a^{s}(z))\}.
\end{eqnarray*}
Similarly, we have
\begin{eqnarray*}
&&\{D_2D_1(x),\a^{k+s}(y),\a^{k+s}(z)\}\\
&=&D'''_2D'''_1(\{x,y,z\})
   -(-1)^{|D_1||x|}\{D_2(\a^{k}(x)), \a^sD'_1(y),\a^{k+s}(z)\}\\
   &&-(-1)^{(|D_1|+|D_2|)|x|}\{\a^{k+s}(x), D'_2 D'_1(y),\a^{k+s}(z)\}\\
   &&-(-1)^{|D_1||x|+|D_2|(|x|+|y|+|D_1|)}\{\a^{k+s}(x), \a^{s}D'_1(y),D''_2(\a^{k}(z))\}\\
   &&-(-1)^{|D_1|(|x|+|y|)}\{D_2(\a^{k}(x)),\a^{k+s}(y),\a^sD''_1(z)\}\\
   &&-(-1)^{|D_1|(|x|+|y|)+|D_2||x|}\{\a^{k+s}(x),D'_2(\a^{k}(y)),\a^sD''_1(z)\}\\
   &&-(-1)^{(|D_1|+|D_2|)(|x|+|y|)}\{\a^{k+s}(x),\a^{k+s}(y),D''_2 D''_1(z)\}\\
   &&-(-1)^{|D_2|(|D_1|+|x|)}\{\a^sD_1(x),D'_2(\a^k(y)),\a^{k+s}(z)\}\\
   &&-(-1)^{|D_2|(|D_1|+|x|+|y|)}\{\a^sD_1(x),\a^{k+s}(y),D''_2(\a^{k}(z))\}.
\end{eqnarray*}
It follows that
\begin{eqnarray*}
&&\{[D_1, D_2](x),\a^{k+s}(y),\a^{k+s}(z)\}\\
&=&\{D_1D_2(x),y,z\}-(-1)^{|D_1||D_2|}\{D_2D_1(x),y,z\}\\
&=&(D'''_1D'''_2-(-1)^{|D_1|+|D_2|}D'''_2D'''_1)\{x,y,z\}\\
&&-(-1)^{|x|(|D_1|+|D_2|)}\{\a^{k+s}(x),(D'_1D'_2-(-1)^{|D_1|+|D_2|}D'_2D'_1)(y),\a^{k+s}(z)\}\\
&&-(-1)^{(|x|+|y|)(|D_1|+|D_2|)}\{\a^{k+s}(x),\a^{k+s}(y),(D''_1D''_2-(-1)^{|D_1|+|D_2|}D''_2D''_1)(z)\}\\
&=&[D'''_1,D'''_2]\{x,y,z\}
   -(-1)^{|x|(|D_1|+|D_2|)}\{\a^{k+s}(x),[D'_1,D'_2](y),\a^{k+s}(z)\}\\
   &&-(-1)^{(|x|+|y|)(|D_1|+|D_2|)}\{\a^{k+s}(x),\a^{k+s}(y),[D''_1,D''_2](z)\},
\end{eqnarray*}
and easily to check that
\begin{eqnarray*}
 &&[[D_1, D_2](x),\a^{k+s}(y)]=[D''_1,D''_2][x,y]-(-1)^{|x|(|D_1|+|D_2|)}[\a^{k+s}(x),[D'_1,D'_2](y)].
\end{eqnarray*}
Obviously, $[D'_1,D'_2],[D''_1,D''_2]$ and $[D'''_1,D'''_2]$ are contained in $End(L)$,
thus $[D_1, D_2]\in GDer_{\a^{k+l}}(L)\subseteq GDer(L)$, that is, $GDer(L)$ is a subalgebra of $End(L)$.
\medskip

(2) For any $D_1\in ZDer_{\a^{k}}(L),D_2\in Der_{\a^s}(L)$ and $x,y,z\in L$, we have
\begin{eqnarray*}
[D_1,D_2]([x,y,z])
=D_1D_2(\{x,y,z\})-(-1)^{|D_1||D_2|}D_2D_1(\{x,y,z\})=0.
\end{eqnarray*}
Also, we have
\begin{eqnarray*}
&&\{[D_1,D_2](x),\a^{k+s}(y),\a^{k+s}(z)\}\\
&=&\{D_1D_2(x),\a^{k+s}(y),\a^{k+s}(z)\}-(-1)^{|D_1||D_2|}\{D_2D_1(x),\a^{k+s}(y),\a^{k+s}(z)\}\\
&=&0-(-1)^{|D_1||D_2|}\{D_2D_1(x),\a^{k+s}(y),\a^{k+s}(z)\}\\
&=&-(-1)^{|D_1||D_2|}(D_2(\{D_1(x),\a^{k}(y),\a^{k}(z)\})-(-1)^{|D_2|(|D_1|+|x|)}\{\a^sD_1(x),D_2(\a^{k}(y)),\a^{k+s}(z)\})\\
&&+(-1)^{|D_1||D_2|}(-1)^{|D_2|(|D_1|+|x|+|y|)}\{\a^sD_1(x),\a^{k+s}(y),D_2(\a^{k}(z))\}\\
&=&0,
\end{eqnarray*}
and easily to check that
\begin{eqnarray*}
[[D_1,D_2](x),\a^{k+s}(y)]=0.
\end{eqnarray*}
It follows that $[D_1,D_2]\in ZDer_{\a^{k+s}}(L)$.
That is, $ZDer(L)$ is an ideal of $Der(L)$. $\hfill \Box$
\medskip

\begin{lemma}
Let $(L, \a)$  be a Hom-Lie-Yamaguti superalgebra,
 then the following statements hold:

(1)~$[Der(L), C(L)]\subseteq C(L)$.

(2)~$[QDer(L), QC(L)]\subseteq QC(L)$.

(3)~$[QC(L), QC(L)]\subseteq QDer(L)$.

(4)~$C(L)\subseteq QDer(L)$.

(5)~$QDer(L)+QC(L)\subseteq GDer(L)$
\end{lemma}
{\bf Proof.}
(1-(4) are easy to prove and omit them, we only check (6). In fact. Let $D_1\in QDer_{\a^{k}}(L), D_2\in QC_{\a^{k}}(L)$. Then there exist $D'_1, D''_1\in End_{\overline{s}}(L)$, for any $x,y, z\in L$, we have
\begin{eqnarray*}
&&[D_1(x), \a^{k}(y)]+(-1)^{s|x|}[\a^{k}(x), D_1(y)]=D_1'([x, y]),\\
&&\{D_1(x),\a^{k}(y),\a^{k}(z)\}+(-1)^{s|x|}\{\a^{k}(x), D_1(y),\a^{k}(z)\}\\
&&\hspace{4cm}+(-1)^{s(|x|+|y|)}\{\a^{k}(x),\a^{k}(y),D_1(z)\}=D_1''(\{x,y,z\}).
\end{eqnarray*}
Thus, for any $x, y, z\in L$, we have
\begin{eqnarray*}
[(D_1+D_2)(x), \a^{k}(y)]&=&[D_1(x), \a^{k}(y)]+[ D_2(x), \a^{k}(y)]\\
&=& D_1'([x, y])-(-1)^{s|x|}[\a^{k}(x), D_1(y)]+(-1)^{s|x|}[ \a^{k}(x),D_2(y)]\\
&=& D_1'([x, y])-(-1)^{s|x|}[\a^{k}(x), (D_1-D_2)(y)],
\end{eqnarray*}
and
\begin{eqnarray*}
&&\{(D_1+D_2)(x), \a^{k}(y), \a^{k}(z)\}\\
&=& \{D_1(x), \a^{k}(y), \a^{k}(z)\}+\{D_2(x), \a^{k}(y), \a^{k}(z)\}\\
&=&D_1''(\{x,y,z\})-(-1)^{s|x|}\{\a^{k}(x), D_1(y),\a^{k}(z)\}-(-1)^{s(|x|+|y|)}\{\a^{k}(x),\a^{k}(y),D_1(z)\}\\
&&+(-1)^{s|x|}\{ \a^{k}(x), D_2(y), \a^{k}(z)\}\\
&=& D_1''(\{x,y,z\})-(-1)^{s|x|}\{\a^{k}(x), (D_1-D_2)(y),\a^{k}(z)\}-(-1)^{s(|x|+|y|)}\{\a^{k}(x),\a^{k}(y),D_1(z)\}.
\end{eqnarray*}
Therefore, $D_1+D_2\in GDer_{\a^{k}}(L)$.
$\hfill \Box$
\medskip

\begin{proposition}
Let $(L, \a)$  be a Hom-Lie-Yamaguti superalgebra, then $QC(L)+[QC(L), QC(L)]$ is a subalgebra of $GDer(L)$.
\end{proposition}
{\bf Proof.}  By Lemma 4.12, (3) and (6), we have
\begin{eqnarray*}
QC(L)+[QC(L), QC(L)]\subseteq GDer(L),
\end{eqnarray*}
and
\begin{eqnarray*}
&&[QC(L)+[QC(L), QC(L)], QC(L)+[QC(L), QC(L)]]\\
&\subseteq&  [QC(L)+GDer(L), QC(L)+[QC(L), QC(L)]]\\
&\subseteq& [QC(L), QC(L)]+[QC(L), [QC(L), QC(L)]]+[QDer(L), QC(L)]\\
&&\hspace{6cm}[QDer(L), [QC(L), QC(L)]].
\end{eqnarray*}
It is easy to verify that $[QDer(L), [QC(L), QC(L)]]\subseteq [QC(L), QC(L)]$ by the Jacobi identity of Hom-Lie algebra. Thus
\begin{eqnarray*}
[QC(L)+[QC(L), QC(L)], QC(L)+[QC(L), QC(L)]]\subseteq QC(L)+[QC(L), QC(L)]\subseteq GDer(L).
\end{eqnarray*}
$\hfill \Box$
\begin{theorem}
Let $(L, \a)$  be a Hom-Lie-Yamaguti superalgebra, $\a$   surjections,
 then $[C(L), QC(L)]\subseteq End(L,Z(L))$.
  Moreover, if $Z(L)=\{0\}$, then  $[C(L), QC(L)]=\{0\}$.
\end{theorem}

{\bf Proof.}
For any $D_1\in C_{\a^{k}}(L), D_2\in QC_{\a^{s}}(L)$ and $x,y,z\in L$, since $\a$ is surjections, there exist $y', z'\in L$ such that $y=\a^{k+s}(y'), z=\a^{k+s}(z')$, we have
\begin{eqnarray*}
[[D_1, D_2](x), y]&=&[[D_1, D_2](x), \a^{k+s}(y')]\\
&=& [D_1D_2(x), \a^{k+s}(y')]-(-1)^{|D_1||D_2|}[D_2D_1(x), \a^{k+s}(y')]\\
&=& D_1([D_2(x), \a^{s}(y')])-(-1)^{|D_1||D_2|}[\a^{s}D_1(x), D_2\a^{k}(y')]\\
&=& D_1([D_2(x), \a^{s}(y')])-D_1([D_2(x), \a^{s}(y')])\\
&=& 0,
\end{eqnarray*}
and
\begin{eqnarray*}
&&\{[D_1,D_2](x),y,z\}\\
&=&\{D_1D_2(x),\a^{k+s}(y'),\a^{k+s}(z')\}-(-1)^{|D_1||D_2|}\{D_2D_1(x),\a^{k+s}(y'),\a^{k+s}(z')\}\\
&=&D_1(\{D_2(x),\a^{s}(y'),\a^{s}(z')\})-(-1)^{|D_1||D_2|}(-1)^{(|D_1|+|x|)|D_2|}\{\a^sD_1(x), D_2(\a^{k}(y')),\a^{k+s}(z')\}\\
&=&(-1)^{|x||D_2|}D_1(\{\a^{s}(x),D_2(y'),\a^{s}(z')\})-(-1)^{|a||D_2|}D_1(\{\a^{s}(x),D_2(y'),\a^{s}(z')\})\\
&=&0.
\end{eqnarray*}
So $[D_1,D_2](x)\subseteq Z(L)$ and therefore $[C(L), QC(L)]\subseteq End(L,Z(L))$.
 Moreover, if $Z(L)=\{0\}$, then  it is easy to see that $[C(L), QC(L)]=\{0\}$.
$\hfill \Box$

\section{ 1-Parameter formal deformations  of Hom-Lie-Yamaguti superalgebras}
\def\theequation{\arabic{section}. \arabic{equation}}
\setcounter{equation} {0}

 Let $(L, \a)$  be a Hom-Lie-Yamaguti superalgebra over $\mathbb{K}$ and $\mathbb{K}[[t]]$  the power series ring in one variable $t$ with coefficients in $\mathbb{K}$.
Assume that $L[[t]]$ is the set of formal series whose coefficients are elements of the vector space $L$.

\begin{definition}
Let $(L, \a)$  be a Hom-Lie-Yamaguti superalgebra.
 A 1-parameter formal deformations of $L$  is a pair of formal power series
 $(f_{t}, g_t)$ of the form
\begin{eqnarray*}
f_{t}=[\c, \c]+\sum_{i\geq 1}f_{i}t^i,~~~~g_{t}=\{\c, \c, \c\}+\sum_{i\geq 1}g_{i}t^i,
\end{eqnarray*}
where each $f_{i}$ is a $\mathbb{K}$-bilinear map $f_{i}: L\times L\rightarrow L$  (extended to be $\mathbb{K}[[t]]$-bilinear) and each $g_{i}$ is a $\mathbb{K}$-trilinear map $g_{i}: L\times L \times L\rightarrow L$  (extended to be $\mathbb{K}[[t]]$-bilinear) such that  $(L[[t]], f_t, g_t, \a)$ is a  Hom-Lie-Yamaguti superalgebra  over $\mathbb{K}[[t]]$.        Set $f_{0}=[\c, \c]$ and $g_0=\{\c, \c, \c\}$, then $f_t$ and $g_t$ can be written as  $f_{t}=\sum_{i\geq 0}f_{i}t^i, g_{t}=\sum_{i\geq 0}g_{i}t^i$, respectively.
\end{definition}
Since $(L[[t]], f_t, g_t, \a)$ is a  Hom-Lie-Yamaguti superalgebra. Then it  satisfies the  following axioms:
\begin{eqnarray}
&& ~\a\circ f_t(x, y)=f_t(\a(x), \a(y))\\
&& ~\a\circ g_t(x,y,z)=g_t(\a(x),\a(y),\a(z))\\
&&~f_t(x, x)=0\\
&&~g_t(x,x,y)=0\\
&& ~\circlearrowleft_{(x,y,z)}(-1)^{|x||z|}(f_t(f_t(x, y), \a(z))+g_t(x, y, z))=0\\
&& ~\circlearrowleft_{(x,y,z)}(-1)^{|x||z|}g_t(f_t(x, y), \a(z), \a(u))=0\\
&&~g_t(\a(x), \a(y), f_t(u, v))=f_t(g_t(x, y, u), \a^{2}(v))+(-1)^{|u|(|x|+|y|)}f_t(\a^{2}(u),g_t(x, y, v))~~~~~~~\\
&&~ g_t(\a^2(x), \a^{2}(y), g_t(u, v, w))=g_t(g_t(x, y, u), \a^{2}(v), \a^2(w))\nonumber\\
&&\hspace{5cm} +(-1)^{|u|(|x|+|y|)}g_t(\a^{2}(u),  g_t(x, y, v), \a^{2}(w))\nonumber\\
&&\hspace{5cm} +(-1)^{(|u|+|v|)(|x|+|y|)}g_t(\a^{2}(u),  \a^{2}(v), g_t(x, y, w))
\end{eqnarray}
for all $ x, y, z, u, v, w\in L$.

\begin{remark}
Equations $(5.1)-(5.8)$ are equivalent to ($n=0,1,2,\cdots$)
\begin{eqnarray}
&& ~\a\circ f_n(x, y)=f_n(\a(x), \a(y))\\
&& ~\a\circ g_n(x,y,z)=g_n(\a(x),\a(y),\a(z))\\
&&~f_n(x, x)=0\\
&&~g_n(x,x,y)=0\\
&& ~\circlearrowleft_{(x,y,z)}(-1)^{|x||z|}(\sum_{i+j=n}f_i(f_j(x, y), \a(z))+g_n(x, y, z))=0\\
&& ~\circlearrowleft_{(x,y,z)}(-1)^{|x||z|}\sum_{i+j=n}g_i(f_j(x, y), \a(z), \a(u))=0\\
&&~\sum_{i+j=n}g_i(\a(x), \a(y), f_j(u, v))=\sum_{i+j=n}f_i(g_j(x, y, u), \a^{2}(v))\nonumber\\
&&\hspace{5cm}+(-1)^{|u|(|x|+|y|)}f_i(\a^{2}(u),g_j(x, y, v))\\
&&~\sum_{i+j=n} g_i(\a^2(x), \a^{2}(y), g_j(u, v, w))=\sum_{i+j=n}g_i(g_j(x, y, u), \a^{2}(v), \a^2(w))\nonumber\\
&&\hspace{4cm} +(-1)^{|u|(|x|+|y|)}g_i(\a^{2}(u),  g_j(x, y, v), \a^{2}(w))\nonumber\\
&&\hspace{4cm} +(-1)^{(|u|+|v|)(|x|+|y|)}g_i(\a^{2}(u),  \a^{2}(v), g_j(x, y, w))
\end{eqnarray}
for all $ x, y, z, u, v, w\in L$,
respectively. These equations are called the deformation equations of a  Hom-Lie-Yamaguti superalgebra. Eqs.(5.9)-(5.12) imply $(f, g)\in C^2(L, L)\times C^3(L, L)$.
\end{remark}

Let $n=1$ in Eqs. (5.8)-(5.16). Then
\begin{eqnarray*}
&& ~\circlearrowleft_{(x,y,z)}(-1)^{|x||z|}([f_1(x, y), \a(z)]+f_1([x, y], \a(z))+g_1(x, y, z))=0\\
&& ~\circlearrowleft_{(x,y,z)}(-1)^{|x||z|}(\{f_1(x, y), \a(z), \a(u)\}+g_1([x, y], \a(z), \a(u)))=0\\
&&~\{\a(x), \a(y), f_1(u, v)\}+g_1(\a(x), \a(y), [u, v])-[g_1(x, y, u), \a^{2}(v)]\nonumber\\
&&-f_1(\{x,y,u\}, \a^{2}(v))-(-1)^{|u|(|x|+|y|)}(\a^{2}(z), g_1(x, y, z))-(-1)^{|u|(|x|+|y|)}f_1(\a^{2}(z),\{x, y, u\})\\
&& \{\a^2(x), \a^{2}(y), g_1(u, v, w)\}+g_1(\a^2(x), \a^{2}(y), \{u, v, w\})-\{g_1(x, y, u), \a^{2}(v), \a^2(w)\}\nonumber\\
&& -g_1(\{x, y, u\}, \a^{2}(v), \a^2(w))-(-1)^{|u|(|x|+|y|)}\{\a^{2}(u),  g_1(x, y, v), \a^{2}(w)\}\nonumber\\
&&-(-1)^{|u|(|x|+|y|)}g_1(\a^{2}(u),  \{x, y, v\}, \a^{2}(w))-(-1)^{|u|(|x|+|y|)}\{\a^{2}(u),\a^{2}(v), g_1(x, y, w) \} \nonumber\\
&&-(-1)^{(|u|+|v|)(|x|+|y|)}g_1(\a^{2}(u),  \a^{2}(v), \{x, y, w\})-(-1)^{(|u|+|v|)(|x|+|y|)}\{\a^{2}(u),  \a^{2}(v), g_1(x, y, w)\}.
\end{eqnarray*}

 which imply $(\delta^{2}_{I}, \delta^{2}_{II})(f_1,g_1)=(0, 0)$, i.e.,
 \begin{eqnarray*}
 (f_1,g_1)\in Z^2(l, L)\times Z^3(l, L).
 \end{eqnarray*}
And $(f_{1}, g_1)$ is called the infinitesimal  deformation of $(f_{t}, g_t)$.
\medskip

\begin{definition}
Let $(L, \a)$  be a Hom-Lie-Yamaguti superalgebra.
 Two 1-parameter formal deformations $(f_{t}, g_t)$ and $(f'_{t}, g'_t)$  of $L$ are said to be equivalent,  denoted by $(f_{t}, g_t)\sim (f'_{t}, g'_t)$,
 if there exists a formal isomorphism of  $\mathbb{K}[[t]]$-modules
  \begin{eqnarray*}
\phi_{t}(x)=\sum_{i\geq 0}\phi_{i}(x)t^{i}:(L[[t]],f_{t},g_t, \a)\rightarrow (L[[t]],f'_{t},g'(t), \a),
\end{eqnarray*}
where $\phi_{i}:L \rightarrow L$ is a $\mathbb{K}$-linear map (extended to be $\mathbb{K}[[t]]$-linear) such that
 \begin{eqnarray*}
&& \phi_{0}=id_L, ~~~~\phi_{t}\circ \a=\a\circ \phi_{t},\\
&&\phi_{t}\circ f_t(x, y)=f'_t(\phi_{t}(x),\phi_{t}(y)), \phi_{t}\circ g_t(x, y, z)=g'_t(\phi_{t}(x),\phi_{t}(y), \phi_{t}(z)).
\end{eqnarray*}
\end{definition}
In particular, if $(f_1, g_1)=(f_2, g_2)=\cdots =(0, 0),$ then  $(f_t, g_t)=(f_0, g_0)$ is called the null deformation.
If  $(f_t, g_t)\sim (f_0, g_0)$, then $(f_t, g_t)$ is called the trivial deformation.
If every 1-parameter formal deformation $(f_t, g_t)$ is trivial, then $L$ is called an analytically rigid Hom-Lie-Yamaguti superalgebra.

\begin{theorem}
Let $(f_{t}, g_t)$ and $(f'_{t}, g'_t)$ be
 two equivalent 1-parameter formal deformations  of $L$.
Then the infinitesimal deformations $(f_1, g_1)$ and $(f'_1, g'_1)$
belong to the same cohomology class in  $H^{2}(L,L)\times H^{3}(L,L).$
\end{theorem}

{\bf Proof.}
By the assumption that $(f_1, g_1)$ and $(f'_1, g_1')$ are equivalent, there exists  a formal isomorphism $\phi_{t}(x)=\sum_{i\geq 0}\phi_{i}(x)t^{i}$
 of  $\mathbb{K}[[t]]$-modules satisfying
  \begin{eqnarray*}
\sum_{i\geq 0}\phi_i(\sum_{j\geq 0}f_j(x_1,x_2)t^{j})t^{i}
=\sum_{i\geq 0}f'_i(\sum_{k\geq 0}\phi_{k}(x_1)t^{k},\sum_{l\geq 0}\phi_{l}(x_2)t^{l})t^{i},\\
\sum_{i\geq 0}\phi_i(\sum_{j\geq 0}g_j(x_1,x_2, x_3)t^{j})t^{i}
=\sum_{i\geq 0}g'_i(\sum_{k\geq 0}\phi_{k}(x_1)t^{k},\sum_{l\geq 0}\phi_{l}(x_2)t^{l}, \sum_{m\geq 0}\phi_{m}(x_3)t^{m})t^{i},\\
\end{eqnarray*}
for any $x_1,x_2, x_3\in L.$
Comparing with the coefficients of $t^1$ for two sides of the above equation, we have
    \begin{eqnarray*}
    &&f_1(x_1,x_2)+\phi_{1}([x_1,x_2])=f'_1(x_1,x_2)+[\phi_{1}(x_1),x_2]+[x_1,\phi_{1}(x_2)],\\
&&g_1(x_1,x_2, x_3)+\phi_{1}(\{x_1,x_2, x_3\})\\
&=&g'_1(x_1,x_2, x_3)+\{\phi_{1}(x_1),x_2, x_3\}+\{x_1,\phi_{1}(x_2), x_3\}+\{x_1,x_2, \phi_{1}(x_3)\}.
\end{eqnarray*}
It follows that $(f_1-f'_1, g_1-g'_1)=(\delta^{1}_{I}, \delta^{1}_{II})(\phi_1, \phi_1)\in B^{2}(L,L)\times B^{3}(L,L),$  as desired.
The proof is completed.$\hfill \Box$

\begin{theorem}
Let $(L, \a)$  be a Hom-Lie-Yamaguti superalgebra with $H^{2}(L,L)\times H^{3}(L,L)=0$, then $L$ is   analytically rigid.
\end{theorem}

{\bf Proof.}
Let $(f_{t}, g_t)$ be a 1-parameter formal deformation  of $L$. Suppose $f_t=f_0+\sum_{i\geq r}f_it^i$, $g_t=g_0+\sum_{i\geq r}g_it^i$.
Set $n=r$ in Eqs.(5.9)-(5.12), It follows that
\begin{eqnarray*}
(f_r,g_r)\in Z^{2}(L,L)\times Z^{3}(L,L)=B^{2}(L,L)\times B^{3}(L,L).
\end{eqnarray*}
Then there exists $h_r\in C^{1}(L, L)$ such that $(f_r,g_r)=(\delta^1_{I}h_r, \delta^1_{I}h_r)$.

Consider $\phi_t=id_L-h_rt^r,$ then $\phi_t: L\rightarrow L$ is a linear isomorphism and $\phi_t\circ \a=\a\circ \phi_t$. Thus we can define another  1-parameter formal deformation  by $\phi_t^{-1}$ in the form of
\begin{eqnarray*}
f'_t(x, y)= \phi_t^{-1}f_t(\phi_t(x), \phi_t(y)), g'_t(x, y, z)= \phi_t^{-1}g_t(\phi_t(x), \phi_t(y), \phi_t(z)).
\end{eqnarray*}
Set $f'_t=\sum_{i\geq 0}f'_it^i$ and use the fact $\phi_tf'_t(x, y)=f_t(\phi_t(x), \phi_t(y))$, then we have
\begin{eqnarray*}
(id_L-h_rt^r)\sum_{i\geq 0}f'_i(x, y)t^i=(f_0+\sum_{i\geq 0}f_it^i)(x-h_r(x)t^r, y-h_r(y)t^{r}),
\end{eqnarray*}
that is
\begin{eqnarray*}
&&\sum_{i\geq 0}f'_{i}(x, y)t^{i}-\sum_{i\geq 0}h_r\circ f'_{i}(x, y)t^{i+r}\\
&=&f_0(x, y)
   -f_0(h_r(x),y)t^{r}-f_0(x,h_r(y))t^r+f_0(h_r(x),h_r(y))t^{2r}\\
   &&+\sum_{i\geq r}f_{i}(x,y)t^{i}-\sum_{i\geq r}\{f_i(h_r(x),y)-f_i(x, h_r(y))\}t^{i+r}+\sum_{i\geq r}f_i(h_r(x),h_r(y))t^{i+2r}.
\end{eqnarray*}
By the above equation, it follows that
\begin{eqnarray*}
&&f'_0(x, y)=f_0(x, y)=[x, y],\\
&&f'_1(x, y)=f'_1(x, y)=\cdots=f'_{r-1}(x, y)=0,\\
&&f'_r(x, y)-h_r([x, y])=f_r(x, y)-[h_r(x),y]-[x,h_r(y)].
\end{eqnarray*}
Therefore, we deduce
\begin{eqnarray*}
f'_r(x, y)=f_r(x, y)-\delta^1_{I}h_r(x, y)=0.
\end{eqnarray*}
It follows that $f'_t=[\c, \c]+\sum_{i\geq r+1}f'_it^i$. Similarly, we have $g'_t=\{\c, \c, \c\}+\sum_{i\geq r+1}g'_it^i$
 By induction, we have $(f_t, g_t)\sim (f_0, g_0)$, that is, $L$ is   analytically rigid.
 The proof is finished.
$\hfill \Box$

\begin{center}
 {\bf ACKNOWLEDGEMENT}
 \end{center}

  The paper is  supported by the NSF of China (No. 11761017), the Youth Project for Natural Science Foundation of Guizhou provincial department of education (No. KY[2018]155), and the Anhui Provincial Natural Science Foundation (No. 1808085MA14).

\renewcommand{\refname}{REFERENCES}

\end{document}